\documentclass[11pt]{amsart}
\usepackage{texmac}
\operators{Stab,Sel,Gal,Disc}
\calsymbols{c}{C,K,O}
\def\theequation{\arabic{equation}}
\def\vK{v_{\scriptscriptstyle K}}
\def\vF{v_{\scriptscriptstyle F}}
\def\Eth{D}
\def\Type{\par\smallskip\noindent Type }
\let\ge\geqslant

\def\disc{d}
\begin{document}

\title[A remark on Tate's algorithm and Kodaira types]{A remark on Tate's algorithm and\\Kodaira types}

\author{Tim and Vladimir Dokchitser}
\address{Dept of Mathematics, University Walk, Bristol BS8 1TW, United Kingdom}
\email{tim.dokchitser@bristol.ac.uk}
\address{Emmanuel College, Cambridge CB2 3AP, United Kingdom}
\email{v.dokchitser@dpmms.cam.ac.uk}

\begin{abstract}
We remark that Tate's algorithm to determine the minimal model of 
an elliptic curve can be stated in a way that characterises Kodaira types
from the minimum of $v(a_i)/i$. As an application, we deduce the behaviour
of Kodaira types in tame extensions of local fields.
\end{abstract}

\maketitle


Let $R$ be a complete discrete valuation ring with 
perfect residue field,~frac-\\tion field~$K$ and valuation $v$.
If $E/K$ is an elliptic curve in Weierstrass form,
$$
  y^2+a_1xy+a_3y = x^3 +a_2x^2+a_4x+a_6 \qquad\qquad (a_i\in K),
$$
the celebrated algorithm of Tate (\cite{TatA}, \cite{Sil2} \S IV.9) determines 
the minimal model and the local invariants of $E$. 
In this paper we gently tweak the resulting models so that the Kodaira type
can be simply read off from the valuations of the $a_i$:

\begin{theorem}
\label{tate}
An elliptic curve $E/K$ with additive reduction has a minimal
Weierstrass model over $R$ which depends on its Kodaira type as follows:

\medskip\noindent
\begin{center}
\begin{tabular}{|r|c|c|c|c|c|c|c|c|c|}
\hline
&\vphantom{$\int^{X^A}$}\II & \III & \IV & \IZS & \InS{n>0} & \IVS & \IIIS & \IIS \\[2pt]
\hline
$\displaystyle\min\frac{v(a_i)}i=$ &
$\displaystyle\vphantom{\int^{X^A}}\frac 16$ & $\displaystyle\frac 14$ 
  & $\displaystyle\frac 13$ & $\displaystyle\frac 12$ & $\displaystyle\frac12$ 
  & $\displaystyle\frac 23$ & $\displaystyle\frac 34$ & $\displaystyle\frac 56$ 
  \\[9pt]
{\rm\footnotesize extra condition}& & & \tiny$v(b_6)\!=\!2$ & \tiny$v(\disc)\!=\!6$ & 
    \rlap{\raise5.5pt\hbox{\tiny$v(\disc)\!>\!6$}}%
    \raise-4pt\hbox{\tiny$v(a_2^2\hbox{-}\hskip 1pt3a_4)\!=\!2$}
& \tiny$v(b_6)\!=\!4$ & & \\[4pt]
\hline
\end{tabular}
\end{center}

\smallskip\noindent
Here
$$
  b_6\!=\!a_3^2\!+\!4a_6\!=\!\Disc(y^2\!+\!a_3y\!-\!a_6), \qquad
  \disc\!=\!\Disc(x^3\!+\!{a_2}x^2\!+\!{a_4}x\!+\!{a_6}).
$$
Conversely, a Weierstrass model satisfying one of these
conditions is minimal, with the corresponding Kodaira type.
\end{theorem}

There is a refinement for type $\InS{n}$ that recovers $n$ as well:

\begin{proposition}
\label{instar}
An elliptic curve $E/K$ with Kodaira type $\InS{n}$, $n>0$ 
has a minimal model with
$$
  v(a_2)=1, \quad
  v(a_i)\ge \tfrac i2 + \lfloor \tfrac {i-1}2\rfloor \tfrac n2,\quad 
  \ch{v(\disc)=n+6,\>\, v(b_6)\ge n+3&\text{ if $\>2\>|\>n,$}\cr
      v(\disc)\ge n+6,\>\, v(b_6)= n+3&\text{ if $\>2\nmid n.$}\cr}
$$
Conversely, a Weierstrass equation satisfying these conditions 
defines an elliptic curve with Kodaira type $\InS{n}$.
\end{proposition}

As an application, we deduce the behaviour 
of minimal discriminants and Kodaira types in tame extensions of local fields;
our motivation came from Iwasawa theory of elliptic curves, where it is necessary
to control local invariants of elliptic curve in towers of number fields, 
see \cite{megasha}.

%
%
%
%
%

\begin{theorem}
\label{tame}
Let $F/K$ be a tame extension of ramification degree $e$, 
and let $E/K$ be an elliptic curve.
\begin{enumerate}
\item If $E/K$ has Kodaira type \In{n}, then $E/F$ has type \In{en}.
\item If $E/K$ has Kodaira type \InS{n}, then $E/F$ has type \InS{en} 
  if $e$ is odd and type \In{en} if $e$ is even.
\item In all other cases, the type of $E/F$ is determined by 
$$
  \eth_{E/F}\equiv e\,\eth_{E/K}\mod 12,
$$
where $\eth=0,2,3,4,6,8,9,10$ if $E$ has Kodaira type 
\IZ, \II, \III, \IV, \IZS, \IVS, \IIIS, \IIS{} respectively.
\end{enumerate}
The valuations of minimal discriminants for $E/K$ and $E/F$ are related by
$$
  v_F(\Delta_{E/F}) = e\,v_K(\Delta_{E/K}) -
        12 \lfloor \frac{e\,\eth_{E/K}}{12}\rfloor,
$$
where 
$\eth_{E/K}=0$ for $\IZ$ and $\In{n}$, 
$6$ for \InS{n} and is as in {\rm (3)} otherwise.
\end{theorem}

\begin{remark*}
If the residue characteristic is at least 5 and $E/K$ has 
potentially good reduction, 
the fraction in the table in Theorem \ref{tate} is just $\frac{v(\Delta_{E/K})}{12}$, 
and $\eth_{E/K}=v(\Delta_{E/K})$ in Theorem \ref{tame}.
The conclusion of Theorem \ref{tame} is then equivalent to the standard fact 
that $\vF(\Delta_{E/F})<12$.  
The point is that $\eth$ gives the correct replacement for $v(\Delta)$ 
in residue characteristics 2 and 3.
Note, however, that in residue characteristics 2 and 3 neither the Kodaira 
type nor the minimal discriminant behave as in Theorem \ref{tame} in 
wild extensions.
\end{remark*}

\begin{example}
The curve $E: y^2=x^3-2x$ over $K=\Q_2$ has Kodaira type~\III{}
($\min\frac{v(a_i)}i=\frac 14$) and $v(\Delta)=9$.
By Theorem \ref{tame}, over the tame extensions $F_n=\Q_2(\sqrt[5^n]2)$ 
the reduction remains of Type \III, and the valuations 
$$
  \displaystyle
  v_{F_n}(\Delta_{E/F_n}) =
  \displaystyle
  9\cdot 5^n - 12 \lfloor \frac{3\cdot 5^n}{12} \rfloor 
    = 6\cdot 5^n+3 
  \displaystyle
   = 33, 153, 753, \ldots. 
$$
In particular, they are not bounded by 12 (or by anything) as they would be
in residue characteristics $\ge 5$.
Over the wild quartic extensions $\Q_2(\sqrt[4]2)$, $\Q_2(\sqrt[4]{-2})$,
$\Q_2(\zeta_8)$, the Kodaira types of $E$ are \IIIS, \InS{3}{}, \InS{4}, and
valuations of the minimal discriminants are 12, 12, 24, respectively.
So these cannot be recovered just from $E/\Q_2$ and the ramification degree.
\end{example}

In the proofs below we follow the steps of Tate's algorithm, numbered 
as in \cite{TatA} and \cite{Sil2} \S IV.9.

\section{Proof of Theorem \ref{tate}}

Let $\pi$ be a uniformiser of $K$. 

By Steps 1-2 of Tate's algorithm, an elliptic curve with additive reduction
over $K$ has a model with $\pi|a_3,a_4,a_6$, $\pi|b_2=a_1^2+4a_2$.
If $K$ has residue characteristic 2, this means $\pi|a_1$, and, 
shifting $y\to y-\alpha x$ for any $\alpha\in\cO_K$ with 
$\alpha^2\equiv a_2\mod\pi$ we can get $\pi|a_2$ as well.
Similarly, if $K$ has odd residue characteristic,  
the substitution $y\mapsto y-\frac{a_1}{2}x$ 
makes both $a_1$ and $a_2$ divisible by $\pi$.
Now we run Tate's algorithm through this equation, and inspect the model
that comes out of it:

Type \II{} (Step 3): Here $\pi^2\nmid a_6$
and the valuations of the $a_i$ are $\ge\!1$, $\ge\!1$, $\ge\!1$, $\ge\!1$, $=\!1$, 
so $\min\frac{v(a_i)}{i}=\frac16$.

Type \III{} (Step 4): Here $\pi^2|a_6$ and 
$$
  \pi^3\nmid b_8 = a_1^2a_6+4a_2a_6-a_1a_3a_4+a_2a_3^2-a_4^2
    \equiv -a_4^2\mod\pi^3.
$$
So $\pi^2\!\nmid\!a_4$, the
valuations of the $a_i$ are $\ge\!1$, $\ge\!1$, $\ge\!1$, $=\!1$, $\ge\!2$ 
and $\min\frac{v(a_i)}{i}\!=\!\frac14$.

Type \IV{} (Step 5): Here $\pi^2|a_6$, $\pi^3|b_8 \Rightarrow \pi^2| a_4$, 
and $\pi^3\nmid b_6 = a_3^2+4a_6$. 
The valuations of the $a_i$ are $\ge 1, \ge 1, \ge 1, \ge 2, \ge 2$,
and either $v(a_3)=1$ or $v(a_6)=2$ since $\pi^3\nmid a_3^2+4a_6$.
So $\min\frac{v(a_i)}{i}=\frac13$.

Type \IZS{} (Step 6): Here $\pi^2|a_3, \pi^2|a_4, \pi^3|a_6$,
$v(\disc)=6$, the
valuations of the $a_i$ are $\ge 1, \ge 1, \ge 2, \ge 2, \ge 3$,
so $\min\frac{v(a_i)}{i}\ge\frac12$.
Because 
$$
  6 = v(\disc) = \pi^6\,\Disc(x^3+\frac{a_2}{\pi}x^2+\frac{a_4}{\pi^2}x+\frac{a_6}{\pi^3}),
$$ 
at least one of $\frac{a_2}{\pi}$, 
$\frac{a_4}{\pi^2}$ and $\frac{a_6}{\pi^3}$ is a unit, so the minimum is 
exactly $\frac 12$.


Type \InS{n}, $n\ge 1$ (Step 7): Here 
$\pi^2|a_3, \pi^2|a_4, \pi^3|a_6$, 
$v(\disc)>6$ and $\pi^2\nmid a_2$,
so $\min\frac{v(a_i)}{i}=\frac12$, attained for $i=2$.
Moreover, the cubic $x^3+\frac{a_2}{\pi}x^2+\frac{a_4}{\pi^2}x+\frac{a_6}{\pi^3}$ 
has a double root which is not a triple root.
A cubic polynomial $x^3+ax^2+bx+c$ has a triple root if and only if its 
discriminant is 0 and $a^2-3b=0$, and this gives the two extra 
stated conditions\footnote
{If the roots of $x^3+ax^2+bx+c$ are $\alpha,\beta,\gamma$,
then the discriminant condition is equivalent to two of them 
being equal, say $\alpha=\beta$, in which case 
$a^2-3b=(\alpha-\gamma)^2$ measures whether 
it is a triple root.}.

Type \IVS{} (Step 8): Here $\pi^2|a_2, \pi^2|a_3, \pi^3|a_4, \pi^4|a_6$
and $\pi^5\nmid b_6 = a_3^2+4a_6$. 
The valuations of the $a_i$ are $\ge 1, \ge 2, \ge 2, \ge 3, \ge 4$,
and either $v(a_3)=2$ or $v(a_6)=4$ since $\pi^5\nmid a_3^2+4a_6$.
So $\min\frac{v(a_i)}{i}=\frac23$ .
%

Type \IIIS{} (Step 9): Here $\pi^2|a_2, \pi^3|a_3, \pi^3|a_4, \pi^5|a_6$ and $\pi^4\nmid a_4$.
The valuations of the $a_i$ are $\ge 1, \ge 2, \ge 3, =3, \ge 5$, 
so $\min\frac{v(a_i)}{i}=\frac34$.

Type \IIS{} (Step 10): Here $\pi^2|a_2, \pi^3|a_3, \pi^4|a_4, \pi^5|a_6$ 
and $\pi^6\nmid a_6$. The 
valuations of the $a_i$ are $\ge 1, \ge 2, \ge 3, \ge 4, =5$, 
so $\min\frac{v(a_i)}{i}=\frac56$.

Conversely, any model satisfying one of the conditions in the table
is minimal with the right Kodaira type, which is immediate from
the corresponding step of Tate's algorithm. 
(The steps do not change such a model.)

\section{Proof of Proposition \ref{instar}}

From Step 7 of Tate's algorithm it follows
readily that a curve $E/K$ of type \InS{n}{} has a minimal model with 
$$
  v(a_2)=1, \quad
  v(a_i)\ge \tfrac i2 + \lfloor \tfrac {i-1}2\rfloor \tfrac n2,\quad 
  \ch{v(D)=n+4,\>\, v(b_6)\ge n+3&\text{ if $\>2\>|\>n,$}\cr
      v(D)\ge n+4,\>\, v(b_6)= n+3&\text{ if $\>2\nmid n,$}\cr}
$$
where $D=\Disc(a_2 x^2 + a_4 x + a_6)$. 
Because $n\ge 1$,
%
%
%
$$
  \disc=-4a_2^3a_6 + a_2^2a_4^2 - 4a_4^3 - 27a_6^2 + 18 a_2a_4a_6 
      \equiv a_2^2 D \mod \pi^{n+7},
$$
so the conditions on $D$ are equivalent to those on $\disc$ in the proposition.

Conversely, such a model has $v(a_2^2-3a_4)=2$, and so the polynomial
$x^3 + \frac{a_2}{\pi}x^2 + \frac{a_4}{\pi^2}x + \frac{a_6}{\pi^3}$
has a double root, but not a triple root. 
Step 7 of Tate's algorithm shows the model to be minimal of type \InS{n}.

\section{Proof of Theorem \ref{tame}}

(1) If $E$ has good or multiplicative reduction (types \IZ, \In{n>0}), 
the minimal model stays minimal in all extensions, and the reduction stays
good, respectively multiplicative. In the multiplicative case, $-n$ is the 
valuation of the $j$-invariant of $E$, so it gets scaled by $e$ in $F/K$;
cf. \cite{Sil2} \S IV.9, Table~4.1.

%
(2), (3) Fix a uniformiser $\pi$ of $F$. 
Write $\vK, \vF$ for the valuations on $K$ and~$F$, and 
$l$ for the residue characteristic. As $F/K$ is tame, $l\nmid e$.

Assume that $E/K$ has additive reduction, and is in Weierstrass form as in 
Theorem \ref{tate}
(and as in Proposition \ref{instar} for type \InS{n}). 
Then $\min \frac{\vK(a_i)}{i}=\frac{\eth_{E/K}}{12}$
and $\min \frac{\vF(a_i)}{i}=\frac{e\,\eth_{E/K}}{12}$. 
Over $F$ this model can be rescaled 
$\lfloor \frac{e\,\eth_{E/K}}{12}\rfloor$ times with the standard
substitution $y\mapsto\pi^3y, x\mapsto\pi^2x$;
call the new Weierstrass coefficients $A_1,A_2,A_3,A_4,A_6$. So now
$$
  \min_i \frac{\vF(A_i)}{i}
   \in \Bigl\{0,\frac 16,\frac 14,\frac 13,\frac 12,\frac 23,\frac 34, \frac 56\Bigr\}.
$$
Now we proceed to show that the resulting equation 
satisfies the `extra conditons' of Theorem \ref{tate}, 
or, if $\min \frac{\vF(A_i)}{i}=0$, that $E/F$ has good reduction (Type \IZ). 
This implies all the claims in the theorem.

A well-known consequence of the fact
that the tame inertia is cyclic is that tame extensions can be built up from
unramified ones and ramified extensions of prime degree.  
If $F/K$ is unramified ($e=1$), there is nothing to prove. So for 
simplicity we may and will assume that $[F:K]=e=p$ is prime, $p\!\ne\! l$.

We first deal with the cases when $E$ acquires good reduction:

\Type \IV, \IVS, $p=3$, $l\!\ne\! 3$: the 
valuations of the $a_i$ are $>\!\frac13,>\!\frac 23,\ge\!1,>\!\frac 43,\ge\! 2$
for Type \IV{} and $>\!\frac23,>\!\frac 43,\ge\! 2,>\!\frac 83,\ge\! 4$ for Type \IVS.
The valuations of the $A_i$ are therefore $>\!0,>\!0,\ge\! 0,>\!0,\ge\! 0$, so 
the model reduces to $y^2+\alpha y=x^3+\beta$ over the residue field of $F$.
It has discriminant $-27(\alpha^2+4\beta)^2$ which is non-zero, since $l\!\ne\! 3$
and $\alpha^2+4\beta\!\ne\! 0$ from the $b_6$ condition for $E/K$.
So $E/F$ has good reduction.

\Type \IZS, $p\!=\!2$, $l\!\ne\! 2$: in the same manner, the valuations 
of the $A_i$ are $>\!0$,~$\ge\!0$, $>\!0$, $\ge\!0$, $\ge\!0$, the
model reduces to $y^2=x^3\!+\!\alpha x^2\!+\!\beta x\!+\!\gamma$ and this
has non-zero discriminant since $l\!\ne\! 2$ and $v_K(\disc)=6$ for $E/K$.

\smallskip
Now we look at the remaining cases, all entirely similar. 

\Type \IV, \IVS{}, $p\ne 3$:
The extra condition in the table for $E/K$ automatically rescales to give 
the one for $E/F$.

\Type \IZS, $p\ne 2$: 
The condition for $E/K$ rescales to give the one for $E/F$.

\Type \II, $p\!=\!2$, $l\ne 2$: $\vK(a_3)\!\ge\! 1, \vK(a_6)\!=\!1$ gives 
$\vF(A_3)\!\ge\! 2$, $\vF(A_6)\!=\!2$. So $\pi^3\nmid A_3^2+4A_6\!=\!B_6$, 
which is the condition for type \IV.

\Type \IIS, $p\!=\!2$, $l\ne 2$: $\vK(a_3)\!\ge\! 3, \vK(a_6)\!=\!5$ gives (after one rescaling)
$\vF(A_3)\!\ge\! 3, \vF(A_6)\!=\!4$. So $\pi^5\nmid A_3^2+4A_6\!=\!B_6$, 
which is the conditon for~\IVS.

\Type \II, $p\!=\!3$, $l\ne 3$: $\vK(a_2)\!\ge\! 1$, $\vK(a_4)\!\ge\! 1$, 
$\vK(a_6)\!=\!1$ gives $\vF(A_2)\!\ge\! 3$, $\vF(A_4)\!\ge\! 3$, $\vF(A_6)\!=\!3$, so 
$$
  x^3 + \frac{A_2}{\pi}x^2 + \frac{A_4}{\pi^2}x + \frac{A_6}{\pi^3}
    \equiv x^3 + \text{unit}\mod\pi,
$$
which has non-zero discriminant as $l\!\ne\!3$.
So $\vF(\Disc(x^3\!+\!{A_2}x^2\!+\!{A_4}x\!+\!{A_6}))\!=\!6$ as required for type $\IZS$.
Type \IIS, $p=3$ is similar.

\Type \III, $p\!=\!2$, $l\ne 2$: $\vK(a_2)\!\ge\! 1, \vK(a_4)\!=\!1, \vK(a_6)\!\ge\! 2$ gives 
$\vF(A_2)\!\ge\! 2$, $\vF(A_4)\!=\!2$, $\vF(A_6)\!\ge\! 4$, so 
$$
  x^3 + \frac{A_2}{\pi}x^2 + \frac{A_4}{\pi^2}x + \frac{A_6}{\pi^3}
    \equiv x^3 + \text{unit}\cdot x\mod\pi,
$$
which has non-zero discriminant as $l\!\ne\!2$.
So $\vF(\Disc(x^3\!+\!{A_2}x^2\!+\!{A_4}x\!+\!{A_6}))\!=\!6$ as required for type $\IZS$.
Type \IIIS, $p=2$ is similar.

\Type \InS{n}, $p=2$, $l\ne 2$ : 
$E$ has non-integral $j$-invariant (\cite{Sil2}~\S~IV.9, Table 4.1),
and so acquires multiplicative reduction over $F$ (\cite{Sil2}~Thm~V.5.3). 
Comparing the valuations of the $j$-invariants and the discriminants, 
we get that $E/F$ has type \In{2n}, and the $A_i$ define a minimal equation.


\Type \InS{n}, $p\ne 2$: The valuations of $a_1,...,a_6,b_6,\disc$
are $\ge\! 1$, $=\!1$, $\ge\!\frac{n+3}2$, $\ge\!\frac{n+4}2$, 
$\ge\!{n\!+\!3}$, $\ge\!{n\!+\!3}$, $\ge\! n\!+\!6$
with equality for one of the last two (depending on whether $n$ is even or odd).
Over $F$ they become $\ge\! p$, $=\!p$, $\ge\! p\frac{n+3}2$, 
$\ge\! p\frac{n+4}2$, $\ge\! p(n\!+\!3)$, $\ge\! p(n\!+\!3)$, $\ge\! p(n\!+\!6)$.
After rescaling the model $\lfloor\frac{6e}{12}\rfloor=\frac{p-1}2$ times, we 
get that the valuations of $A_1,...,A_6,B_6,\Eth$ for the new model are 
$\ge\!\frac{p+1}2$, $=\!1$, $\ge\!\frac{pn+3}2$, $\ge\!\frac{pn+4}2$, 
$\ge\!{pn\!+\!3}$, $\ge\!{pn\!+\!3}$, $\ge\! pn\!+\!6$,
again with equality for one of the last two. 
In other words, it satisfies the conditions of Proposition \ref{instar} for 
Type \InS{pn}.

\begin{acknowledgements}
The first author is supported by a 
Royal Society University Research Fellowship.
We would like to thank Matthias Schuett for carefully
reading the manuscript.
\end{acknowledgements}

\end{document}